\documentclass[letterpaper, 10 pt, conference]{ieeeconf}  % Comment this line out if you need a4paper

\IEEEoverridecommandlockouts                              % This command is only needed if 
\usepackage{graphics} % for pdf, bitmapped graphics files
\graphicspath{{./figures/}}
\usepackage{epsfig} % for postscript graphics files
\usepackage{amsmath} % assumes amsmath package installed
\usepackage{amssymb}  % assumes amsmath package installed
\usepackage{amsfonts}
\usepackage{color}
\usepackage{algorithm}
\usepackage{algpseudocode}
\usepackage{hyperref}
\usepackage{soul}

\DeclareMathOperator*{\argmin}{arg\,min}

\DeclareMathOperator*{\Minimize}{Minimize:}

\DeclareMathOperator*{\SubjectTo}{Subject\phantom{a}to:}
\DeclareMathOperator*{\supp}{\textsf{supp}}

\makeatletter
\def\BState{\State\hskip-\ALG@thistlm}
\makeatother

\newcommand{\norm}[1]{\left\lVert#1\right\rVert} 
\newcommand{\abs}[1]{\left\lvert#1\right\rvert}

\newtheorem{lemma}{\bf{Lemma}}
\newtheorem{definition}{\bf{Definition}}
\newtheorem{proposition}{\bf{Proposition}}
\newtheorem{theorem}{\bf{Theorem}}
\newtheorem{corollary}{\bf{Corollary}}
\newtheorem{remark}{Remark}

\usepackage{mathtools}
\DeclarePairedDelimiter{\ceil}{\lceil}{\rceil}

\title{\LARGE \bf
Robust Resilient Signal Reconstruction under Adversarial Attacks
}

\author{Yu Zheng$^{1}$~\IEEEmembership{GSIEEE}, Olugbenga Moses Anubi$^{1}$~\IEEEmembership{SMIEEE}, Lalit Mestha$^{2}$~\IEEEmembership{FIEEE},
Hema Achanta$^{3}$~\IEEEmembership{MIEEE}% <-this % stops a space
\thanks{*This work was not supported by any organization}% <-this % stops a space
\thanks{$^{1}$Yu Zheng and Olugbenga Moses Anubi are with the Department of Electrical and Computer Engineering, Florida State University, FL, USA.}%
\thanks{$^{2}$Lalit Mestha is with Genetic Innovations Inc., Honolulu, HI, USA (work performed while at GE Global Research)}%
\thanks{$^{3}$Hema Achanta is with GE Global Research, Niskayuna, NY, USA.}
\thanks{O. Anubi is the corresponding author, {\tt\small oanubi@fsu.edu}}.
}

\begin{document}

\maketitle
\thispagestyle{empty}
\pagestyle{empty}

%%%%%%%%%%%%%%%%%%%%%%%%%%%%%%%%%%%%%%%%%%%%%%%%%%%%%%%%%%%%%%%%%%%%%%%%%%%%%%%%
\begin{abstract}
We consider the problem of signal reconstruction for a system under sparse signal corruption by a malicious agent. The reconstruction problem follows the standard error coding problem that has been studied extensively in the literature. We include a new challenge of robust estimation of the attack support. The problem is then cast as a constrained optimization problem merging promising techniques in the area of deep learning and estimation theory. A pruning algorithm is developed to reduce the ``false positive" uncertainty of data-driven attack localization results, thereby improving the probability of correct signal reconstruction. Sufficient conditions for the correct reconstruction and the associated reconstruction error bounds are obtained for both exact and inexact attack support estimation. Moreover, a simulation of a water distribution system is presented to validate the proposed techniques.
\end{abstract}

%%%%%%%%%%%%%%%%%%%%%%%%%%%%%%%%%%%%%%%%%%%%%%%%%%%%%%%%%%%%%%%%%%%%%%%%%%%%%%%%
\section{INTRODUCTION}\label{Introduction}
% **MOTIVATION FOR NEUTRALIZATION *****\\
The majority of Industrial systems and critical infrastructures are Cyber-physical Systems (CPS), in that they consist of an interplay between physical components (sensors, controllers, and actuators) and digital components (computational algorithms, software systems, human-machine interfaces) via communication networks \cite{lee2015cyber}. This opens up a portal that makes them prime targets of cyber malicious activities \cite{pasqualetti2013attack}. The resilient signal reconstruction is a filtering problem for removing undesirable effects created by malicious intent as in adversarial attacks or other similar unbounded phenomena on some of the system’s monitoring nodes \cite{fawzi2014secure}. For physical systems, if signal reconstruction is performed jointly on all the affected signals, then it can improve the resiliency of critical infrastructure to cyber-activities and fault-induced anomalies to allow continued, safe operation \cite{mestha2017cyber}. Conventional resilient estimation designs, such as event-trigger Luenberger-like observer \cite{shoukry2015event, lu2019switched}, constrained sensor fusion \cite{nakahira2018attack,chen2018resilient}, $\ell_1$ decoders \cite{pajic2017design, Zheng2022ResilientPrior}, generally admit an assumption that no more than half of measurements are compromised. Recent results show this assumption could be relaxed if the prior knowledge is incorporated into the estimation scheme \cite{anubi2019resilient, shinohara2019resilient, khazraei2022attack}. 

One of the good sources of prior information is active attack detection and localization (AADL) which could produce a guess of the location of attacks. AADL is a promising technique where people consider how a defense policy could intelligently distinguish the attacks from the nominal signals. In \cite{liu2020online}, the authors leveraged a random input as a watermark and identified attacks by examining the outputted watermark. In \cite{anubi2019enhanced}, the authors added a data-driven authentication model to a power system by training a Gaussian regression model based on the power market data. Moreover, different data-driven algorithms or hybrid physics-data-driven algorithms are also employed to explore the underlying difference between attacked signals and nominal signals. In \cite{esmalifalak2014detecting}, the authors implemented a distributed support vector machine to identify stealth attacks in projected feature space. Unsupervised/semi-supervised approaches were also widely used to obtain good generalization to unseen attacks, such as generative adversarial network (GAN) \cite{ferdowsi2019generative}, GAN-based multi-layer perceptron classifier \cite{Zheng2021algorithm}, deep autoencoding Gaussian mixture model-based detector \cite{ahmed2020challenges}. 

In general, AADL is seen as a procedure to estimate the attack support. Research in compressed sensing has demonstrated significant improvement to sparse recovery performance by incorporating support information. For example, relaxing the sparsity assumption and reducing the recovery error \cite{friedlander2011recovering}. However, due to the black-box nature of data-driven learning models, AADL's performance relies on the quality of the training dataset and the choice of pre-defined hyperparameters \cite{ahmed2020challenges}. Thus, AADL could only produce probabilistic conclusions that would contain ``false positive" (FP) cases. Consequently, the FP cases in the estimated attack support should be quantified and utilized in the signal reconstruction scheme. In this paper, we quantify the FP uncertainty of AADL using a Bernoulli distribution model and develop a symmetrical pruning operation to compensate for the FP uncertainty. Next, a robust resilient signal reconstruction method is proposed using the pruning algorithm.

% This paper is an updated version of our arxiv paper with the same title \cite{anubi2018robust}. The arxiv version of this paper was never published or submitted to any conference or journal although the idea presented in the paper has already been cited. Beyond the arxiv version of the paper, we complete the theoretical developments of the pruning algorithm for robust resilient signal reconstruction with inexact support knowledge. We also add a simulation to give an example of how to use the proposed techniques. Moreover, the introduction section was improved to clarify the motivation and contributions of this paper.

% *******PAPER ORGANIZATION GO HERE********\\
The rest of the paper is organized as follows. Necessary notations are clarified in Section \ref{notation}. In Section \ref{Prelims}, a measurement model is presented with the basic reconstruction problem. In Section \ref{Reconstruction}, some achievable error bounds are proved for the reconstruction problem where the exact support of the attack vector is assumed to be provided apriori by some support estimator. In Section \ref{pruning}, the exact support knowledge assumption is removed and a more realistic scenario is considered: the estimated attack support is assumed to satisfy a Bernoulli distribution. A pruning algorithm is given to compensate for the support uncertainty and new error bounds are obtained based on the pruning algorithm. A simulation of a water distribution system was studied in Section \ref{sim}. Finally, concluding remarks and future directions are highlighted in Section \ref{conclusion}.

\section{NOTATION}\label{notation}
We use the same notations of real numbers, real vectors, and real matrices as in \cite{zheng2021attack}. Other necessary notations are clarified here for the subsequent development in the rest of this paper. The support of a vector $\mathbf{x}$ is given by $\textsf{supp}(\mathbf{x})=\{i|\mathbf{x}_i\neq0\}$. Given a set $\mathcal{S} \subseteq \{1,\cdots, n\}$, $X_{\mathcal{S}} \in {\mathbb R}^{|\mathcal{S}| \times n}$ is the matrix containing the rows of the matrix $X \in \mathbb{R}^{m\times n}$ indexed by $\mathcal{S}$. Also, $\mathcal{S}_1 \backslash \mathcal{S}_2 \triangleq \mathcal{S}_1 \cap \mathcal{S}_2^c$ denotes set difference. $\mathcal{S}^c$ denotes the complement of the set $\mathcal{S}$ in a universal set. An indicator vector $\mathbf{q} \in \{0, 1\}^{n}$ of $\mathcal{S} \subseteq \{1,\dots, n\}$ is given as:
 \begin{equation} \label{equ:support}
     \mathbf{q}_{i}=\left\{
     \begin{array}{lr}
     0 &\text{if}\hspace{0.2cm} i \in \mathcal{S} \\
     1 &\text{otherwise}.
     \end{array}
     \right.
\end{equation}
We use $\textsf{argsort} \downarrow (\mathbf{x})$ to represent a function returning a vector of indices $i$ of vector elements $\mathbf{x}_i$ in the descending order of the magnitude of $\mathbf{x}_i$. The symbols $\&$, $*$, and $\odot$ denote the logical AND operator, the convolution multiplication, and the \textcolor{blue}{point-wise} multiplication respectively. $\mathbb{E}(x)$ denotes the expected value of a random variable $x$. A Bernoulli distributed variable $x$ with probability $\textsf{Pr}\{x=1\}=p$ is denoted by $x \sim \mathcal{B}(1,p)$. To measure the similarity of two binary vectors, a positive predictive value (PPV) is used in this paper and defined as follows.
\begin{definition}\label{Def:PPV_def}
\textbf{(Precision, Positive Prediction Value, $\textsf{PPV}$ \cite{fawcett2006introduction})}
Given two binary vectors $\hat{\mathbf{x}}, \mathbf{x} \in\{0,1\}^{n}$, in that order, $\textsf{PPV} $ is the ratio of the entries of $\mathbf{x}$ correctly estimated in $\hat{\mathbf{x}}$. It is given by
\begin{equation}\label{equ:PPV_def}
 \textsf{PPV}(\hat{\mathbf{x}},\mathbf{x}) = \frac{TP}{TP+FP} = \frac{\|\mathbf{x} \odot \hat{\mathbf{x}} \|_{0}}{\max\{1,\|\hat{\mathbf{x}}\|_{0}\}} ,
\end{equation}
\textcolor{blue}{where $TP$ is true positive defined as an outcome where the model correctly predicts the positive class and $FP$ is true negative defined as an outcome where the model correctly predicts the negative class}
\end{definition}
\iffalse
\noindent We shall use $\textsf{PPV}$ to indicate the function $\textsf{PPV}(\hat{\mathbf{x}},\mathbf{x})$ whence the meaning of $\hat{\mathbf{x}},\mathbf{x}$ are clear from the content. Last, a lemma is given
\begin{lemma}[\cite{fernandez2010closed}]  \label{Lem:PBD_CDF}
The sum of $N$ independent Bernoulli distributed variables $\epsilon_i \sim \mathcal{B}(1,\mathbf{p}_i)$ satisfies the below pdf
\begin{equation}\label{equ:PMF}
      Pr\left\{\sum_{i=1}^{N} \epsilon_{i} = k-1 \right\}=\mathbf{r}_k, \hspace{2mm} k=1,\cdots,N+1,
\end{equation} 
where $ \mathbf{r} = \beta \cdot \begin{bmatrix}-\mathbf{s}_1\\1 \end{bmatrix} * \begin{bmatrix}-\mathbf{s}_2\\1 \end{bmatrix} * \cdots * \begin{bmatrix}-\mathbf{s}_m\\1 \end{bmatrix} \in \mathbb{R}^{N+1}$ with $ \beta =\displaystyle\prod_{i=1}^{N} \mathbf{p}_{i}, \hspace{2mm} \mathbf{s}_i = -\frac{1-\mathbf{p}_{i}}{\mathbf{p}_{i}}$.
\end{lemma}
\fi

% $\mathbb{R}$ denotes the space of real numbers, $\mathbb{R}^n$ denotes the space of real vectors of length $n$, and $\mathbb{R}^{m\times n}$ real matrices of $m$ rows and $n$ columns. $\mathbb{R}_+$ denotes positive real numbers. $X^\top$ denotes the transpose of the matrix $X$. We use normal-face lower-case letters ($x\in\mathbb{R}$) to represent real scalars, use bold-face lower-case letters ($\mathbf{x}\in\mathbb{R}^n$) to represent vectors, and use normal-face upper-case letters ($X\in\mathbb{R}^{n\times m}$) to represent matrices.

\section{PRELIMINARIES}\label{Prelims}

Due to the attacker's limited resources, they are only capable of compromising a sparse set of sensors simultaneously at any time \cite{pasqualetti2013attack}. By assuming the sparsity of the attack vector, the adversarial signal reconstruction problem could be formulated as an error correction problem \cite{pajic2017design}. Given a coding matrix $F\in\mathbb{R}^{n\times m}$ ($n\ll m$), and a measurement vector $\mathbf{y}\in\mathbb{R}^m$, the attack recovery problem is to recover a sparse vector \textcolor{blue}{$\mathbf{e}\in \mathbb{R}^m$,} $\left\|\mathbf{e}\right\|_{0}<m$ subject to $\mathbf{y}=F\mathbf{e}$. This problem has been studied in compressive sensing:
\begin{align}\label{eqn:comp_sens}
\Minimize\limits_{\mathbf{e}}{\left\|\mathbf{e}\right\|_0}\hspace{2mm}\SubjectTo\hspace{2mm}\mathbf{y}=F\mathbf{e}.
\end{align}
Hayden et. al \cite{hayden2016sparse} found a condition of achieving unique reconstruction of  $\left\|\mathbf{e}\right\|_0\le q$ by \eqref{eqn:comp_sens}: all sub-columns of $2q$ columns of $F$ are full-rank. This condition is also called ``$2q$-observerability" in resilient estimation literature \cite{chong2015observability, shoukry2017secure}. Although the program in \eqref{eqn:comp_sens} is solved as is in some cases~\cite{pajic2017design}, it is in fact an NP-hard problem due to its nonconvexity and index counting objective. Consequently, it is usually relaxed to a convex optimization program:
\begin{align}\label{eqn:comp_sens_L1}
\Minimize\limits_{\mathbf{e}}{\left\|\mathbf{e}\right\|_{1}}\hspace{2mm}\SubjectTo\hspace{2mm}\mathbf{y}=F\mathbf{e}.
\end{align}
The program \eqref{eqn:comp_sens_L1} produces the equivalent solution if a restricted isometric property (RIP) of matrix $F$ holds~\cite{donoho2003optimally,elad2002generalized}, given by
\begin{align}\label{equ:RIP}
(1-\delta_k)\norm{\mathbf{v}}^2\le\norm{F({\mathcal{T})}\mathbf{v}}^2\le(1+\delta_k)\norm{\mathbf{v}}^2
\end{align}
for all subsets $\mathcal{T}$ with $\abs{\mathcal{T}}\le k$ and vector $\mathbf{v}\in\mathbb{R}^{\abs{\mathcal{T}}}$, and $\delta_k$ is called $k$-restricted isometry constant defined as the minimum quantity under which RIP \eqref{equ:RIP} holds. If $F$ satisfies RIP, then its every sub-columns of cardinality less than $k$ are approximately orthonormal. Furthermore, it has been proved that if 
\begin{align*}
\delta_k + \delta_{2k} + \delta_{3k}<1,
\end{align*}
\textcolor{blue}{where $\delta_{2k}$ and $\delta_{3k}$ are defined similarly to \eqref{equ:RIP} with $|\mathcal{T}|\leq 2k$ and $|\mathcal{T}|\leq 3k$ respectively}. Then any sparse signal $\mathbf{e}$, with $\abs{\textsf{supp}(\mathbf{e})}\le k$, could be reconstructed by \eqref{eqn:comp_sens}.

Now, suppose there exists an oracle that estimates the support $\mathcal{T}$ in advance, then the sparse vector $\mathbf{e}$ can be estimated to an accuracy of $\frac{2\epsilon}{1-\delta_{\abs{\mathcal{T}}}}$\footnote{$\norm{\hat{\mathbf{y}}-F({\mathcal{T})}\hat{\mathbf{e}}}^2\leq \norm{\hat{\mathbf{y}}-F({\mathcal{T})}\mathbf{e}}^2 =\norm{\hat{\mathbf{y}}-\mathbf{y}}^2 \leq \epsilon$. Thus, $\norm{\hat{\mathbf{y}}-\hat{\mathbf{y}}+F({\mathcal{T})}(\mathbf{e}-\hat{\mathbf{e}})}^2 \leq \epsilon \Rightarrow \norm{F({\mathcal{T})}(\mathbf{e}-\hat{\mathbf{e}})}^2\leq 2\epsilon \Rightarrow (1-\delta_{\abs{\mathcal{T}}})\norm{\mathbf{e}-\hat{\mathbf{e}}}^2\leq 2\epsilon$} by the least square estimator:
\begin{align*}
\hat{\mathbf{e}} = \argmin\limits_{\mathbf{e}}\left\{\norm{\hat{\mathbf{y}}-F({\mathcal{T})}\mathbf{e}}^2\right\},
\end{align*}
where $\hat{\mathbf{y}}$ is the measured signal with $\|\hat{\mathbf{y}}-\mathbf{y}\|^2\leq \epsilon$, and $\epsilon$ is the associated measurement error. Of course, if the measurement is error-free, then the estimation is exact. This work aims to investigate the least-square-type estimator when such an oracle is available subject to both oracle and measurement uncertainties. Such oracles are termed \emph{localization oracles} that AADL usually generates.

Consider the linear model:
\begin{align}\label{eqn:linear_model}
\mathbf{y} = C\mathbf{x} + \mathbf{e} + \mathbf{v},
\end{align}
where $\mathbf{y},\mathbf{e},\mathbf{v} \in \mathbb{R}^m$ are vectors of measurements, attack/corruption due to an adversarial agent, and error term due to measurement noise/model uncertainty respectively. The matrix $C\in\mathbb{R}^{m\times n}$ is a mapping from some internal state ($\subseteq\mathbb{R}^n$) to the output space ($\subseteq\mathbb{R}^m$). The following assumptions are made with respect to the model above \cite{pajic2017design,nakahira2015dynamic,zheng2021Resilient}:
\begin{description}
	\item[$\bullet$ Redundancy:] {\hspace{1.3cm}Measurements contain redundant information in that $m>n$ }
	\item[$\bullet$ Bounded Noise:] {\hspace{1.7cm}There exists a known $\varepsilon>0$ such that $\norm{\mathbf{v}}\le\varepsilon$}
	\item[$\bullet$ Sparse Corruption:] {\hspace{2cm} $\textsf{supp}(\mathbf{e})\ll m$}
	\item[$\bullet$ Attack-Noise Orthogonality:] {\hspace{3.5cm} Without loss of generality, $\mathbf{e}^\top\mathbf{v}=0$}
\end{description}
 
Consequently, the reconstruction problem is given by:

\begin{align}\label{eqn:recon}
\begin{array}{rl}
\Minimize  & \norm{\mathbf{e}}_0 + \norm{\mathbf{v}}_{2}\\
\SubjectTo &\\
           &\mathbf{y}=C\mathbf{x} + \mathbf{e} + \mathbf{v}\\
					 &\mathbf{e}^\top\mathbf{v} = 0
\end{array}
\end{align}
which is, in general, a very challenging problem to solve due to the index minimization objective and the degeneracy introduced by the complementarity constraint $\mathbf{e}^\top\mathbf{v} = 0$. However, if there exists a localization oracle that provides the support $\mathcal{T} = \textsf{supp}(\mathbf{e})$ apriori, then the reconstruction problem reduces to the unconstrained problem:
\begin{align}
\Minimize\norm{\mathbf{y}_\mathcal{T}-C_\mathcal{T}\mathbf{x}}_{2}.
\end{align}
Of course, there are obvious conditions under which the solution to the above optimization problem provides no guarantee of reconstructing the actual signal. In what follows, the reconstruction error bounds are studied in more detail under different conditions.

\section{Reconstruction with Exact Support Knowledge}\label{Reconstruction}
In this section, we examine some bounds on the reconstruction error when the attack support is known exactly. Although the exact knowledge assumption is not pragmatic, it does give us a lower bound and a benchmark for the cases where the support is not known exactly. The next result gives the performance of a least-square reconstruction from partial information.
\begin{theorem}[Least Square Reconstruction]\label{thm:LS}
Given the linear model
\begin{align}\label{eqn:global_LS}
\mathbf{y} = C\mathbf{x} + \boldsymbol{\nu},
\end{align}
where $\mathbf{y}\in\mathbb{R}^m$ is a vector of measurements, $\mathbf{x}\in\mathbb{R}^n, n\le m$ is a vector of internal states (or features), $C\in\mathbb{R}^{m\times n}$, and $\boldsymbol{\nu}$ is the model error with the associated error bound $\|\boldsymbol{\nu}\|\le\varepsilon$ for a known constant $\varepsilon>0$.
\\\\\noindent
Consider any partial measurement $\mathbf{y}_1\in\mathbb{R}^{m_1},m_1>n$ satisfying
\begin{align}\label{eqn:global_LS1}
\mathbf{y}_1 = C_1\mathbf{x}^* + \boldsymbol{\nu}_1,
\end{align}
where $C_1\in\mathbb{R}^{m_1\times n}$ is a matrix of the corresponding rows of $C$ and $\boldsymbol{\nu}_1$ is the associated model error, the vector $\mathbf{x}^*\in \mathbb{R}^{n}$ is the unknown actual internal state associated with the complete measurement set  as in~\eqref{eqn:global_LS}.
\\\\\noindent
The least-square estimator
\begin{align}\label{equ:LSE_exact}
\hat{\mathbf{x}}=\argmin{\left\{\frac{1}{2}\left\|\mathbf{y}_1-C_1\mathbf{x}\right\|^2\right\}},
\end{align}
of $\mathbf{x}^*$, satisfies the error bound
\begin{align}
\left\|\hat{\mathbf{x}}-\mathbf{x}^*\right\|\le\frac{2}{\sigma_1}\varepsilon,
\end{align}
where $\sigma_1$ is the minimal singular value of $C_1$.
 
\end{theorem}
 
\begin{proof}
Based on \eqref{equ:LSE_exact}, it follows that
\begin{align}
\left\|\mathbf{y}_1-C_1\hat{\mathbf{x}}\right\|^2\le\left\|\mathbf{y}_1-C_1\mathbf{x}^*\right\|^2.
\end{align}
After using~\eqref{eqn:global_LS1}, the above inequality can be simplified and expanded as follows:
\begin{align}\nonumber
&\left\|C_1\mathbf{x}^*-C_1\hat{\mathbf{x}}+\boldsymbol{\nu}_1\right\|^2\le\left\|\boldsymbol{\nu}_1\right\|^2\\
\Rightarrow&\left\|C_1\widetilde{\mathbf{x}}\right\|^2\le 2\boldsymbol{\nu}_1^\top\left(C_1\widetilde{\mathbf{x}}\right),
\end{align}
where $\widetilde{\mathbf{x}} = \hat{\mathbf{x}} - \mathbf{x}^*$. After using Young's Inequality, for some $\delta>0$, the above inequality yields
\begin{align}\nonumber
&\left\|C_1\widetilde{\mathbf{x}}\right\|^2\le \delta\left\|\boldsymbol{\nu}_1\right\|^2 + \frac{1}{\delta}\left\|C_1\widetilde{\mathbf{x}}\right\|^2, \\\nonumber
\Rightarrow&\left(1-\frac{1}{\delta}\right)\left\|C_1\widetilde{\mathbf{x}}\right\|^2\le\delta\left\|\boldsymbol{\nu}_1\right\|^2\\
\Rightarrow&\left\|C_1\widetilde{\mathbf{x}}\right\|^2\le\frac{\delta^2}{\delta-1}\left\|\boldsymbol{\nu}_1\right\|^2.
\end{align}
From which we conclude that (after setting $\delta=2$)
\begin{align}
\left\|\widetilde{\mathbf{x}}\right\|\le\frac{2}{\sigma_1}\left\|\boldsymbol{\nu}\right\|\le\frac{2}{\sigma_1}\varepsilon
\end{align}
\end{proof}

\begin{remark}[Rank-deficiency and RIP]
Necessarily \textcolor{blue}{$\abs{\mathcal{T}^c}\geq n$}, otherwise the reconstruction error $\norm{\hat{\mathbf{x}}-\mathbf{x}^*}$ is unbounded. Consequently, one can conclude that: $\norm{\hat{\mathbf{x}}-\mathbf{x}^*}\le\frac{2}{\delta_n}\varepsilon$, where $\delta_n$ is the $n$-restricted isometry constant of $C^\top$.
\end{remark}
Although it was shown in~\cite{candes2005decoding} that random matrices satisfy the RIP condition with overwhelming probability, certifying such property is still an NP-hard problem~\cite{bandeira2013certifying}. In order to guarantee bounded reconstruction error for the cases where there is a potential loss of row-rank after selection due to the localization oracle, we investigate the use of a special constraint in the reconstruction optimization problem. 

\begin{corollary}\label{cor:constrained_LS}
\textbf{(Constrained Least Square Reconstruction)} Let $\mathcal{X}\subset\mathbb{R}^n$ be a set characterized by $\norm{\mathbf{x}_1-\mathbf{x}_2}\le \delta$ for all $\mathbf{x}_1,\mathbf{x}_2\in\mathcal{X}$ and some $\delta>0$. Consider the constrained least-square estimator:
\begin{align}
\hat{\mathbf{x}}=\argmin\limits_{\mathbf{x}\in\mathcal{X}}{\left\{\frac{1}{2}\left\|\mathbf{y}_1-C_1\mathbf{x}\right\|^2\right\}}.
\end{align}
If $\mathbf{x}^*\in\mathcal{X}$, then the reconstruction error can be bounded as:
\begin{align}
\norm{\hat{\mathbf{x}}-\mathbf{x}^*}\le2\min\left\{\frac{\delta}{2},\frac{\varepsilon}{\delta_n}\right\}.
\end{align}
\end{corollary}
\begin{proof}
Using the optimality of $\hat{\mathbf{x}}$ and following similar argument as in the proof of Theorem~\ref{thm:LS}, it is shown that $\norm{\mathbf{x}-\mathbf{x}^*}\le2\frac{\varepsilon}{\delta_n}$. Next, using the feasibility of both $\mathbf{x}$ and $\mathbf{x}^*$, it follows that $\norm{\mathbf{x}-\mathbf{x}^*}\le\delta$. Thus, $\norm{\hat{\mathbf{x}}-\mathbf{x}^*}\le\min\left\{\delta,2\frac{\varepsilon}{\delta_n}\right\}$.
\end{proof}

\begin{remark}
Although Corollary~\ref{cor:constrained_LS} provides a guaranteed bound on the reconstruction error, it introduces another challenge of finding the set $\mathcal{X}$ that contains the unknown vector $\mathbf{x}^*$. Fortunately, there is a host of supervised and unsupervised machine learning models and algorithms that can be used to find such a set from historical data, together with some exogenous supporting measurement. In such cases, the bound is guaranteed with a probability depending on the \textcolor{blue}{receiver of characteristics (ROC)} of the underlying machine learning model. Interested readers are directed to the reference \cite{mestha2017cyber} where supervised learning with support vector machines was used to generated a local approximation for $\mathcal{X}$ via a quadratic approximation of the boundary score function. 
\end{remark}

\section{Reconstruction with INEXACT Support Knowledge}\label{pruning}
Even though the constrained least square reconstruction described in the previous section provides a guaranteed bound on the reconstruction error, it is impractical due to the exact knowledge assumption on the support estimation. Exact support knowledge alludes to a perfect localization oracle which is not possible, or extremely challenging at best, from a practical standpoint. Thus, it is imperative to understand the effect of the imperfection of the localization oracle on the reconstruction error bound. The goal of this section is to re-examine the constrained least square reconstruction with uncertain localization information. 

Let $\mathcal{T}=\supp(\mathbf{e})$ be the unknown support with a corresponding indicator vector $\mathbf{q}$ defined as in \eqref{equ:support}. Suppose there exists an AADL giving estimated support $\hat{\mathcal{T}}$, with $\hat{\mathbf{q}}$ similarly defined. \textcolor{blue}{Since the localization of attacks is a binary classification procedure, the $FP$ uncertainty of AADL can be described by the following poisson binomial model}:
\begin{align}\label{equ:uncertainty_model}
\mathbf{q}_i = \varepsilon_i\hat{\mathbf{q}}_i + (1-\varepsilon_i)(1-\hat{\mathbf{q}}_i)
\end{align}
 where $\varepsilon_i \sim \mathcal{B}(1,\mathbf{p}_i)$ \textcolor{blue}{is the agreement between the actual support $\mathcal{T}^c$ and the AADL's output $\hat{\mathcal{T}}^c$}, and the probability $\mathbf{p}_i  \in (0,1)$ can be calculated as true rate $\mathbf{p}_i = \mathbb{E}(\varepsilon_i)$ by ROC. Next, the following lemma gives a way to calculate the precision of the underlying AADL based on the above uncertainty model.

\begin{lemma}\label{Lem:PPV_math}
If the estimated support prior $\hat{\mathcal{T}}$ satisfies the uncertainty model in \eqref{equ:uncertainty_model}, its precision could be calculated as
\begin{equation}\label{equ:PPV_math}
    \textsf{PPV} = \frac{1}{|\hat{\mathcal{T}}^c|}\sum_{i \in \hat{\mathcal{T}}^c}\epsilon_i.
\end{equation}
\end{lemma}
\begin{proof}
Based on the definition of indicator vector in \eqref{equ:support}, it holds that $\hat{\mathbf{q}}_i\hat{\mathbf{q}}_i = \hat{\mathbf{q}}_i$. Thus, multiplying $\hat{\mathbf{q}}_i$ on both sides of \eqref{equ:uncertainty_model} yields
$$\begin{aligned}
    \mathbf{q}_i \hat{\mathbf{q}}_i &= \left(\varepsilon_i\hat{\mathbf{q}}_i + (1-\varepsilon_i)(1-\hat{\mathbf{q}}_i)\right)\hat{\mathbf{q}}_i \\
     & = 2\varepsilon_i\hat{\mathbf{q}}_i\hat{\mathbf{q}}_i + \hat{\mathbf{q}}_i - \varepsilon_i\hat{\mathbf{q}}_i - \hat{\mathbf{q}}_i\hat{\mathbf{q}}_i = \varepsilon_i\hat{\mathbf{q}}_i
\end{aligned} 
$$
This implies that
$$\begin{aligned}
\textsf{PPV} =  \frac{\|\mathbf{q} \odot \hat{\mathbf{q}} \|_{0}}{\|\hat{\mathbf{q}}\|_{0}} = \frac{\sum_{i=1}^{m} \mathbf{q}_i \hat{\mathbf{q}}_i}{\sum_{i=1}^{m} \hat{\mathbf{q}}_i}=\frac{1}{|\hat{\mathcal{T}}^c|}\sum_{i=1}^{m} \epsilon_i \hat{\mathbf{q}}_i \\
= \frac{1}{|\hat{\mathcal{T}}^c|}\sum_{i \in \hat{\mathcal{T}}^c}\epsilon_i.
\end{aligned}$$
\end{proof} 

\noindent If a binary classification algorithm cannot outperform a random flip of a fair coin, it is often not seen as a successful algorithm. Thus, the following result gives a condition for AADL outperforming a random flip of a fair coin based on the uncertainty model in \eqref{equ:uncertainty_model}.
\begin{proposition}\label{Proposition:EPPV_0.5}
The underlying AADL outperforms a random flip of a fair coin if and only if 
\begin{equation}\label{equ:assum_condition_pruning}
    \sum_{i=1}^{m} \mathbf{p}_i > m p_A
\end{equation}
where $p_A \in (0,1)$ is the expected percentage of attacked nodes. Furthermore, if the maximum percentage of attacked nodes is $P_A$, then \eqref{equ:assum_condition_pruning} is only a sufficient condition.
\end{proposition}
\begin{proof}
Expanding \eqref{equ:uncertainty_model} yields 
\begin{equation*}
\mathbf{q}_i = 2\epsilon_i \hat{\mathbf{q}}_i +1 -\hat{\mathbf{q}}_i-\epsilon_i = 1-\hat{\mathbf{q}}_i-\epsilon_i+2\mathbf{q}_i \hat{\mathbf{q}}_i.
\end{equation*} 
This implies that
$ \epsilon_i -1 +\mathbf{q}_i = 2(\mathbf{q}_i \hat{\mathbf{q}}_i-\frac{1}{2}\hat{\mathbf{q}}_i)
$. Summing over $i=1,\cdots, m$ and taking the expected value of both sides yield
\begin{equation*} \sum_{i=1}^{m} \mathbf{p}_i -m +E[\|\mathbf{q}\|_{\ell_0}] = 2 E\left[\|\mathbf{q} \odot \hat{\mathbf{q}} \|_{\ell_0}-\frac{1}{2}\|\hat{\mathbf{q}}\|_{\ell_0}\right].
\end{equation*}
Thus,
$ E\left[\|\mathbf{q} \odot \hat{\mathbf{q}} \|_{\ell_0}-\frac{1}{2}\|\hat{\mathbf{q}}\|_{\ell_0}\right] >0$ if and only if
$\sum_{i=1}^{m} \mathbf{p}_i >m -E[\|\mathbf{q}\|_{\ell_0}] \triangleq mp_A$.
Moreover, if $m -E[\|\mathbf{q}\|_{\ell_0}] \leq m p_A$, it is straightforward to see that $\sum_{i=1}^{m} \mathbf{p}_i > m p_A$ is only sufficient for $E\left[\|\mathbf{q} \odot \hat{\mathbf{q}} \|_{\ell_0}+\frac{1}{2}\|\hat{\mathbf{q}}\|_{\ell_0}\right] >0$.
\end{proof}

\noindent Next, we shall present a statistics-based method to reduce the FP uncertainty of the estimated support prior.
\begin{definition}[Pruning, Pruning algorithm, $\textsf{PPV}_{\eta}$]\label{Def: pruning}
A pruning algorithm is a procedure returning a subset support prior $\hat{\mathcal{T}}_{\eta}^c \subset \{1,\cdots,m\}$ of $\hat{\mathcal{T}}^c$ satisfying 
\begin{equation}\label{equ:origin_pruning_def}
    \hat{\mathcal{T}}_{\eta}^c \subseteq \hat{\mathcal{T}}^c.
\end{equation}  
And the corresponding precision of the pruned support prior can be calculated as
\begin{equation}  \label{equ:PPV_eta_math}
    \textsf{PPV}_{\eta} = \frac{\sum_{i\in \hat{\mathcal{T}}_{\eta}^c}\epsilon_i}{|\hat{\mathcal{T}}_{\eta}^c|}.
\end{equation}
\end{definition}

\begin{proposition}\label{proposition:tradeoff}
Given $\gamma_0 >0$, then
\begin{equation*} \textsf{Pr}\left\{\textsf{PPV}_{\eta} - \gamma_0\textsf{PPV} > 0 \right\}>0 \hspace{1mm}\text{if and only if}\hspace{1mm} \gamma_0 |\hat{\mathcal{T}}_{\eta}^c| < |\hat{\mathcal{T}}^c|.
\end{equation*}
\end{proposition}
\begin{proof}
\begin{equation*} \begin{aligned}
\textsf{PPV}_{\eta} - \gamma_0\textsf{PPV} = \frac{1}{|\hat{\mathcal{T}}_{\eta}^c|}\sum_{i\in \hat{\mathcal{T}}_{\eta}^c}\epsilon_i - \frac{\gamma_0}{|\hat{\mathcal{T}}^c|}\sum_{j \in \hat{\mathcal{T}}^c}\epsilon_j \\
=  \left(\frac{1}{|\hat{\mathcal{T}}_{\eta}^c|}-\frac{\gamma_0}{|\hat{\mathcal{T}}^c|}\right)\sum_{i\in \hat{\mathcal{T}}_{\eta}^c}\epsilon_i -\frac{\gamma_0}{|\hat{\mathcal{T}}^c|}\sum_{j \in \hat{\mathcal{T}}^c\backslash \hat{\mathcal{T}}_{\eta}^c}\epsilon_j .
\end{aligned}
\end{equation*}
Thus, a necessary condition for $\textsf{PPV}_{\eta} - \gamma_0\textsf{PPV} > 0$ for some $\varepsilon_i$ is $\frac{1}{|\hat{\mathcal{T}}_{\eta}^c|}-\frac{\gamma_0}{|\hat{\mathcal{T}}^c|} >0$. Moreover, if $\frac{1}{|\hat{\mathcal{T}}_{\eta}^c|}-\frac{\gamma_0}{|\hat{\mathcal{T}}^c|} >0$, one can find some $\varepsilon_i$ for which $\textsf{PPV}_{\eta} - \gamma_0\textsf{PPV} > 0$.

\end{proof}
\begin{remark}
This proposition shows that the smaller the size of the pruned safe set, the bigger the achievable $\textsf{PPV}$ improvement. However, as will be shown later, the smaller the probability of attaining that improvement. This demonstrates a tradeoff between the possibility and probability of $\textsf{PPV}$ improvement, depending on the pruning aggressiveness.
\end{remark}

\noindent Next, the pruning algorithm in \textcolor{blue}{Algorithm~\ref{Algorithm: Pruning}} refines the estimated safe set $\hat{\mathcal{T}}^c$.

\begin{algorithm}

\uppercase{\romannumeral1}.\hspace{2mm}\textbf{}
Obtain the maximum quantity $l_{\eta}$ of safe channels that are localized by $\hat{\mathcal{T}}^c$ correctly with a probability of at least $\eta \in (0,1)$:
\begin{equation}\label{equ:trust_number}
\begin{aligned}
    l_{\eta} = \max \left\{ |\mathcal{I}| \biggl\lvert \prod\limits_{i\in \mathcal{I}}\mathbf{p}_i \geq \eta, \hspace{2mm}\mathcal{I}\in \hat{\mathcal{T}}^c\right\}.
\end{aligned}
\end{equation}

\uppercase{\romannumeral2}.\hspace{2mm}\textbf{}
Use the current localization prior $\hat{\mathbf{q}}$ and the AADL's historical performance $\mathbf{p}$ to extract the $l_{\eta}$ safest nodes as follows.
\begin{equation}    \label{equ:pruning_alg}
    \hat{\mathcal{T}}^{c}_{\eta} =  \big\{\textsf{argsort} \downarrow (\mathbf{p} \odot \hat{\mathbf{q}})\big\}_{1}^{l_{\eta}}.
\end{equation}
where, $\{\cdot\}_{1}^{l_{\eta}}$ is an index extraction from the first elements to $l_{\eta}$ elements.

\caption{A Robust Pruning Algorithm}
\label{Algorithm: Pruning}
\end{algorithm}

\begin{theorem}\label{Thm:PPV_eta_1}
Suppose there exists an AADL generating estimated support prior $\hat{\mathcal{T}}^c$ satisfying \eqref{equ:uncertainty_model}, through the \textcolor{blue}{\emph{Algorithm \ref{Algorithm: Pruning}} }, the precision of the pruned support prior $\hat{\mathcal{T}}_{\eta}^c$ satisfies
\begin{equation*} \textsf{Pr}\left\{\textsf{PPV}_{\eta} =1 \right\} \geq \eta.
\end{equation*}
\end{theorem}

\begin{proof}
According to \eqref{equ:PPV_eta_math},
\begin{equation*}\begin{aligned}
\textsf{Pr}\left\{\textsf{PPV}_{\eta} =1 \right\} =\textsf{Pr}\left\{\sum_{i\in \hat{\mathcal{T}}_{\eta}^c}\epsilon_i = |\hat{\mathcal{T}}_{\eta}^c|\right\} = \prod\limits_{i\in \hat{\mathcal{T}}_{\eta}^c}\mathbf{p}_i
\end{aligned}
\end{equation*}
Based on \eqref{equ:trust_number}, and \eqref{equ:pruning_alg}, we have $|\hat{\mathcal{T}}_{\eta}^c| = |\mathcal{I}| =l_{\eta}$. Also, since the $\mathbf{p}_i$'s in $\hat{\mathcal{T}}_{\eta}^c$ are the $l_{\eta}$ largest probabilities in $\hat{\mathcal{T}}^c$, it follows that $\prod\limits_{i\in \hat{\mathcal{T}}_{\eta}^c}\mathbf{p}_i \geq \prod\limits_{i\in \mathcal{I}}\mathbf{p}_i \geq \eta
$. Thus $\textsf{Pr}\left\{\textsf{PPV}_{\eta} =1 \right\} \geq \eta$.
\end{proof}

Finally, the reconstruction is given by:
\begin{align}\label{equ:LSE_pruning}
\hat{\mathbf{x}}_{\eta}=\argmin\limits_{\mathbf{x}\in\mathcal{X}}{\left\{\frac{1}{2}\left\|\mathbf{y}_{\hat{\mathcal{T}}_{\eta}^c}-C_{\hat{\mathcal{T}}_{\eta}^c}\mathbf{x}\right\|^2\right\}}.
\end{align}

\begin{theorem}
\textbf{(Least Square Reconstruction with Prior Pruning)} Consider the linear measurement model given in \eqref{eqn:linear_model}. Suppose there exists an AADL that gives an estimate, $\hat{\mathcal{T}}$, of $\supp(\mathbf{e})$ with uncertainty described by \eqref{equ:uncertainty_model}. Given a parameter $\eta \in (0, 1]$ with corresponding quantity $l_{\eta}$ given by \eqref{equ:trust_number}, let $\hat{\mathcal{T}}_{\eta}$ be a new support with the indicator $\hat{\mathbf{q}}_{\eta}$ defined by \eqref{equ:pruning_alg}. If $l_\eta-|\operatorname{supp}(\mathbf{e})| \geq n$, then, with a probability of at least $\eta$, the least-square estimator \eqref{equ:LSE_pruning} satisfies the error bound
\begin{equation}
    \left\|\hat{\mathbf{x}}_\eta-\mathbf{x}^*\right\| \leq 2 \min \left\{\frac{\delta}{2}, \frac{\varepsilon}{1-\delta_n}\right\},
\end{equation}
\end{theorem}

\begin{proof}
For the sake of convenience, we state at the beginning of this proof that all claims made next hold with a probability of at least $\eta$. According to Theorem \ref{Thm:PPV_eta_1}, $\left\|\mathbf{y}_{\hat{\mathcal{T}}_{\eta}^c}-C_{\hat{\mathcal{T}}_{\eta}^c}\mathbf{x}\right\| \leq \epsilon$. Using the optimality of the estimator and following the same procedure as in the proof of Theorem \ref{thm:LS} yields
$$\left\|\hat{\mathbf{x}}_\eta-\mathbf{x}^*\right\| \leq 2 \min \left\{\frac{\delta}{2}, \frac{\varepsilon}{\sigma_\eta}\right\},$$
where $\sigma_\eta$ is the smallest singular value of $C_{\hat{\mathcal{T}}_{\eta}^c}$. Next, since $|\hat{\mathcal{T}}_{\eta}| \leq |\hat{\mathcal{T}}|$ and at least $l_\eta$ nodes are correctly localized by $\hat{\mathcal{T}}$, the following holds
\begin{equation}
\begin{aligned}
l_\eta-\left|\hat{\mathcal{T}}_\eta^c\right| &=l_\eta-\left(m-\left|\hat{\mathcal{T}}_\eta\right|\right) \\ & \leq l_\eta-(m-|\hat{\mathcal{T}}|) =l_\eta-\left|\hat{\mathcal{T}}^c\right|  \\
& \leq \supp(\mathbf{e}).
\end{aligned}    
\end{equation}
Using $l_\eta-|\supp(\mathbf{e})| \geq n$, it follows that
$$l_\eta-\left|\hat{\mathcal{T}}_\eta^c\right| \leq l_\eta-n$$
which implies that $\left|\hat{\mathcal{T}}_\eta^c\right| \geq n$. Thus $\sigma_\eta$ can be lower bounded as $\sigma_\eta\geq (1-\delta_n)$, from which the result follows.
\end{proof}

\section{Simulation}\label{sim}
In this section, we use some numerical simulations to validate the theoretical developments in the proceeding sections. 

\begin{figure}
    \centering
    \includegraphics[scale=0.5]{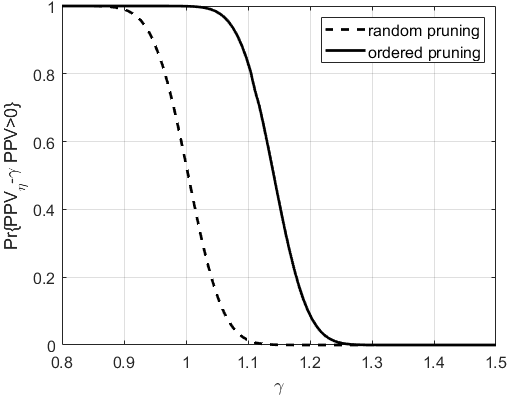}
    \caption{The tradeoff lines between the probability of precision improvement and the possibility $(\gamma)$ of improving precision through pruning operations}
    \label{fig:tradeoff}
\end{figure}

Firstly, to visualize the tradeoff between the possibility and probability shown in Proposition~\ref{proposition:tradeoff}, we perform a random pruning operation and an ordered pruning operation following the Definition~\ref{Def: pruning}. The random pruning operation is to randomly choose a subset of $\hat{\mathcal{T}}^c$ as $\hat{\mathcal{T}}_{\eta}^c$, while the ordered pruning operation is to choose a subset of $\hat{\mathcal{T}}^c$ with the biggest confidences $Pr\{\varepsilon = 1\}$. As shown in figure \ref{fig:tradeoff}, the tradeoff lines between \textcolor{blue}{the probability and possibility $(\gamma)$ of precision improvement} and the possibility $(\gamma)$ are similar s-curves for both pruning operations. Moreover, figure \ref{fig:tradeoff} also shows that ordered pruning outperforms random pruning by considering AADL's confidence. 

\begin{figure}
    \centering
    \includegraphics[scale=0.5]{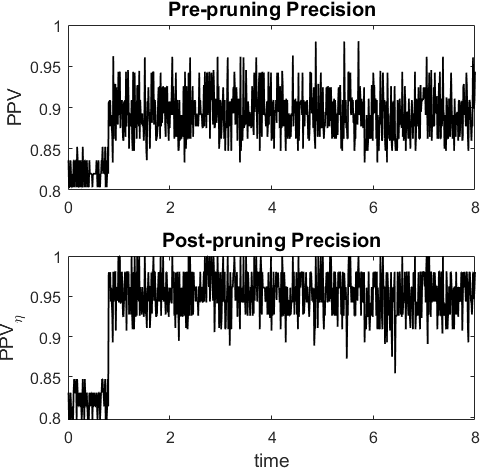}
    \caption{The localization precision of AADL's support prior (pre-pruning precision) and the pruned support prior (post-pruning precision).}
    \label{fig:precision_plot}
\end{figure}

Next, to show the positive correlation between reducing the FP uncertainty in AADL's localization output and better reconstruction performance, we compare the estimation errors of the least-square estimator with the estimated support prior $\hat{\mathcal{T}}^c$ and with the pruned support prior  $\hat{\mathcal{T}}_{\eta}^c$. The comparison is performed on a time-invariant linear model of a water distribution system containing $10$ tanks: 
\begin{equation}\label{equ:water_level_adjustment}
\begin{aligned}
   \mathbf{x}_{i+1} &= A\mathbf{x}_i + B\mathbf{u}_i,\\
    \mathbf{y}_i &= C\mathbf{x}_i  + \mathbf{e}_i + \mathbf{v}_i.
\end{aligned}
\end{equation}
where $\mathbf{x} \in \mathbb{R}^{10}$ is the state vector representing liquid level of water tanks, $\mathbf{u} \in \mathbb{R}^{10}$ is the control signal of magnetic valves, $\mathbf{y}, \mathbf{v} \in \mathbb{R}^{61}$ are the measurements along the pipelines and sensor noise, and $\mathbf{e} \in \mathbf{R}^{61}$ is a modeling sparse vector of attack signals injected through those IoT sensors. The goal of this system is to maintain the water levels in a stable range. The detailed description and control design could be found in \cite{Zheng2021algorithm}. A moving-horizon least-square estimator (LSE) is used to estimate the states:
\begin{align}\label{equ:state_estimate}
    \hat{\mathbf{x}}_i = A^TH_0^\dagger\mathbf{y}_{I} + \left(F-A^TH_0^\dagger H_1\right)\mathbf{u}_{I-1},
\end{align}
where, $\mathbf{y}_I =  [\mathbf{y}_{i-T+1}^{\top} \hspace{1mm} \cdots \hspace{1mm} \mathbf{y}_{i}^{\top} ]^{\top}$, $\mathbf{u}_{I-1} =  [\mathbf{u}_{i-T+1}^{\top} \hspace{1mm}\cdots \hspace{1mm}\mathbf{u}_{i-1}^{\top} ]^{\top}$, $F = [A^{T-1}B \hspace{1mm} A^{T-2}B \hspace{1mm}\hdots \hspace{1mm} B]$,
$$ H_0  =  \begin{bmatrix}CA\\ C A^2 \\ \vdots \\ CA^{T}\end{bmatrix}, H_1  = \begin{bmatrix}CB & 0 & \cdots& 0 \\ CAB & CB& \cdots&  0\\ \vdots & \vdots & \ddots &  \vdots\\ CA^{T-1}B & CA^{T-2}B & \cdots & CB\end{bmatrix}.
$$
The feasible false data injection attack (FDIA) is designed as shown in \cite{zheng2021attack}. In this simulation, approximately $20\%$ ($12$ of $61$) sensors are attacked and the attack location is randomized at each time instance. We trained a suite of Gaussian process regressors (GPR) for each measurement node, as shown in \cite{anubi2019enhanced}. The goal of the AADL algorithm is to detect the attacks if the measurements could not be explained by the associated GPR with high likelihood, as shown in \cite{Zheng2022ResilientPrior}. After detection, the proposed pruning algorithm is used to update the estimated support prior. Figure ~\ref{fig:precision_plot} shows the precision of AADL's localization prior (pre-pruning precision) and the pruned support prior's precision (post-pruning precision). It is seen that the localization precision is improved through the proposed pruning algorithm. 

\begin{figure}
    \centering
    \includegraphics[scale=0.5]{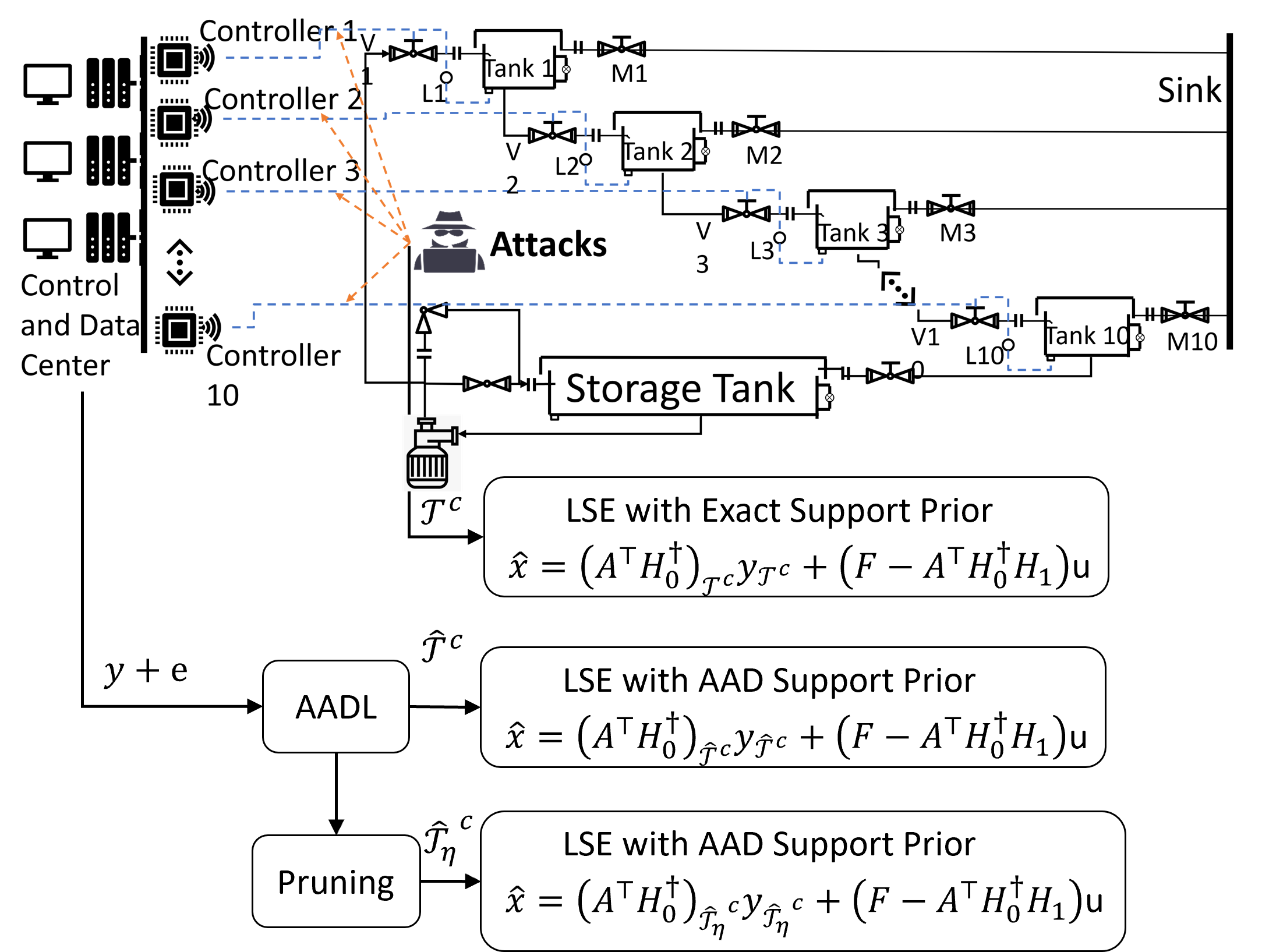}
    \caption{The schematic diagram of the simulation based on a water distribution system}
    \label{fig:sim_scheme}
\end{figure}

\begin{figure}[h!]
    \centering
    \includegraphics[scale=0.5]{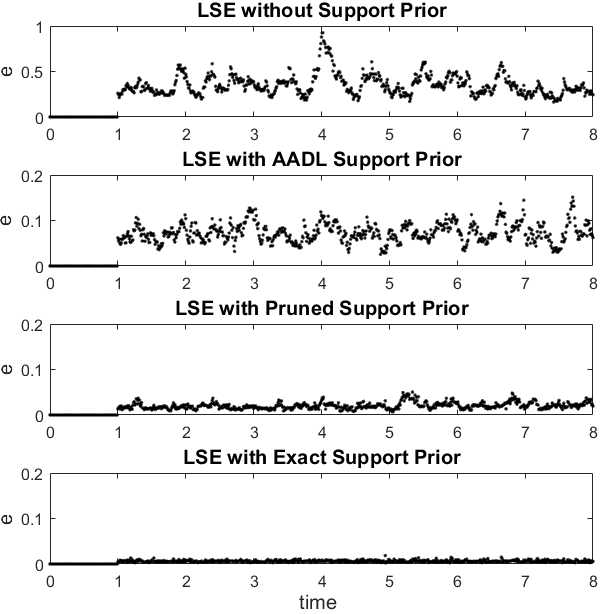}
    \caption{The estimation errors of three different robust resilient estimation schemes}
    \label{fig:error_plot}
\end{figure}

Next, we compared the estimation errors of four signal reconstruction schemes, as shown in Figure~\ref{fig:sim_scheme}\footnote{We open-sourced the simulation at \url{https://github.com/ZYblend/Robust_Resilient_Observer}}. The first one performs an impractical estimation scheme using LSE \eqref{equ:state_estimate} with exact prior knowledge of attack locations. The second one uses the estimated support prior from AADL to clean the measurement nodes used in LSE. The third one uses the proposed pruned support prior. Then, we checked the error of the estimated states from the nominal estimated state without attacks $e = \|\hat{\mathbf{x}}_{\text{with attack}}-\hat{\mathbf{x}}_{\text{without attack}}\|_2$. The estimation errors are shown in Figure~\ref{fig:error_plot}. The mean square errors (MSE) from top to bottom are $0.3085$, $0.0639$, $0.0176$, and $0.0054$. It is seen that the support prior's precision is bigger, then the estimation error is smaller. Without support prior, LSE does not have attack resiliency and thus has the biggest error. With AADL, the estimated support prior is around $90\%$ as shown in \textcolor{blue}{Figure}~\ref{fig:precision_plot}. Then the estimation error is reduced significantly. Although the space of FP uncertainty is already small, the pruning algorithm could still improve the precision of support prior to about $95\%$. Thus, the estimation error is further reduced. The last one shows the perfect estimation performance if the support prior is exact, but the exact support prior could not be known by a defender in practice.

\section{CONCLUSIONS}\label{conclusion}
In this paper, we discuss robust resilient signal reconstruction with attack support prior. A practical scenario, that is the localization prior contains FP uncertainty, is considered. A pruning algorithm is proposed to evaluate and quantify the attack detection and localization algorithm's performance, then improve the localization precision. Although we only present the robust resilient signal reconstruction based on the least-square estimator, a similar estimation scheme with other $2$ norm-based estimators is also expected, such as the Luenberger observer. Moreover, a dynamic update law of observer gain with respect to the pruning algorithm is still an open question.

%\addtolength{\textheight}{-12cm}   % This command serves to balance the column lengths
                                  % on the last page of the document manually. It shortens
                                  % the textheight of the last page by a suitable amount.
                                  % This command does not take effect until the next page
                                  % so it should come on the page before the last. Make
                                  % sure that you do not shorten the textheight too much.

%%%%%%%%%%%%%%%%%%%%%%%%%%%%%%%%%%%%%%%%%%%%%%%%%%%%%%%%%%%%%%%%%%%%%%%%%%%%%%%%

%%%%%%%%%%%%%%%%%%%%%%%%%%%%%%%%%%%%%%%%%%%%%%%%%%%%%%%%%%%%%%%%%%%%%%%%%%%%%%%%

%%%%%%%%%%%%%%%%%%%%%%%%%%%%%%%%%%%%%%%%%%%%%%%%%%%%%%%%%%%%%%%%%%%%%%%%%%%%%%%%
\iffalse
\section*{APPENDIX}

Appendixes should appear before the acknowledgment.

\section*{ACKNOWLEDGMENT}

The preferred spelling of the word ÒacknowledgmentÓ in America is without an ÒeÓ after the ÒgÓ. Avoid the stilted expression, ÒOne of us (R. B. G.) thanks . . .Ó  Instead, try ÒR. B. G. thanksÓ. Put sponsor acknowledgments in the unnumbered footnote on the first page.
\fi

%%%%%%%%%%%%%%%%%%%%%%%%%%%%%%%%%%%%%%%%%%%%%%%%%%%%%%%%%%%%%%%%%%%%%%%%%%%%%%%%
\bibliographystyle{IEEEtran}
% argument is your BibTeX string definitions and bibliography database(s)
\bibliography{references}
\end{document}